\documentclass[12pt]{amsart}

\usepackage{amssymb,amsfonts,amsmath,amsthm,epsfig,tikz,bm,color}

\makeatletter
\def\@seccntformat#1{\csname the#1\endcsname. }
\def\@biblabel#1{#1.}

\makeatother
  
\def\Rm#1{\lowercase\expandafter{\romannumeral#1}}

   \def\Ga{\Gamma}
  \def\D{\Delta} 
  
\def\a{\alpha} \def\b{\beta} \def\g{\gamma}  \def\e{\varepsilon}

 \def\bbZ{{\Bbb Z}}  
   \def\ZZ{\mathbb Z}

  \def\nd{\mathrel{\bigm|\kern-.7em/}}

 \def\Aut{\hbox{\rm Aut}}

\def\Cay{\hbox{\rm Cay}} 
  
\def\qed{\hfill $\Box$} 
  \def\Arc{\hbox{\rm Arc}} 
 \def\bbZ{{\Bbb Z}}  
\def\Hol{\hbox{\rm Hol}}
\def\D2n{\hbox{\rm D$_{2n}$}}

\newtheorem{thm}{Theorem}[section]
\newtheorem{cor}[thm]{Corollary}
\newtheorem{conj}[thm]{Conjecture}
\newtheorem{lem}[thm]{Lemma}

\makeatletter \@addtoreset{equation}{section}

\newtheorem{pro}[thm]{Proposition}

\def\pf{\noindent {\it Proof.\ }}
\def\qed{\ifmmode\square\else\nolinebreak\hfill
$\Box$\fi\par\vskip12pt}

\begin{document}

\title[ Normal Cayley digraphs of dihedral groups with CI-property]
{Normal Cayley digraphs of dihedral groups with CI-property}%
\author{Jin-Hua Xie, $^*$Yan-Quan Feng, Jin-Xin Zhou}%
\address{Department of Mathematics, Beijing Jiaotong University, Beijing, 100044, China\mbox{\hskip 3cm}}
\email{J.-H. Xie, jinhuaxie@bjtu.edu.cn; Y.-Q. Feng, yqfeng@bjtu.edu.cn;\linebreak \mbox{}\hskip 3.75cm J.-X. Zhou, jxzhou@bjtu.edu.cn} %
\thanks{$^*$Corresponding author}

\begin{abstract}

A Cayley (di)graph $\Cay(G,S)$ of a group $G$ with respect to $S$ is said to be normal if the right regular representation of $G$ is normal in the automorphism group of $\Cay(G,S)$, and is called a CI-(di)graph if there is $\a\in \Aut(G)$ such that $S^\a=T$, whenever $\Cay(G,S)\cong \Cay(G,T)$ for a Cayley (di)graph $\Cay(G,T)$. A finite group $G$ is called a DCI-group or a NDCI-group if all Cayley digraphs or normal Cayley digraphs of $G$ are CI-digraphs, and is called a CI-group or a NCI-group if all Cayley graphs or normal Cayley graphs of $G$ are CI-graphs, respectively.

Motivated by a conjecture proposed by \'Ad\'am in 1967, CI-groups and DCI-groups have been actively studied during the last fifty years by many researchers in algebraic graph theory. It takes about thirty years to obtain the classification of cyclic CI-groups and DCI-groups, and recently, the first two authors, among others, classified cyclic NCI-groups and NDCI-groups. Even though there are many partial results on dihedral CI-groups and DCI-groups, their classification is still elusive.
In this paper, we prove that a dihedral group of order $2n$ is a NCI-group or a NDCI-group if and only if $n=2,4$ or $n$ is odd. As a direct consequence, we have that if a dihedral group $\D2n$ of order $2n$ is a DCI-group then $n=2$ or $n$ is odd-square-free, and that if $\D2n$ is a CI-group then $n=2,9$ or $n$ is odd-square-free, throwing some new light on  classification of dihedral CI-groups and DCI-groups.
\end{abstract}

\maketitle
\qquad {\textsc k}{\scriptsize \textsc {eywords.}} {\footnotesize  Dihedral group, CI-group, DCI-group, NCI-group, NDCI-group.}

\section{Introduction}

Graphs and digraphs considered in this paper are finite and simple, and groups are finite. For a (di)graph  $\Gamma$, we use $V(\Gamma)$, $E(\Gamma)$, $\Arc(\Gamma)$ and $\Aut(\Gamma)$ to denote the vertex set, edge set, arc set, and automorphism group of $\Gamma$ respectively, where an arc means a directed edge in a digraph and an order pair of adjacent vertices in a graph.

Let $G$ be a group and $S$ a subset of $G$ with $1\not\in S$. A digraph with vertex set $G$ and arc set $\{(g,sg)\ |\ g\in G, s\in S\}$, denoted by $\Cay(G,S)$, is called the {\em Cayley digraph} of $G$ with respect to $S$. If $S$ is inverse-closed, that is, $S=S^{-1}:=\{x^{-1}\mid x\in S\}$, then for  two adjacent vertices $u$ and $v$ in $\Cay(G,S)$, both $(u,v)$ and $(v,u)$ are arcs, and in this case, we view $\Cay(G,S)$ as a graph by identifying the two arcs with one edge $\{u,v\}$.

Two Cayley (di)graphs $\Cay(G,S)$ and $\Cay(G,T)$ are called {\em Cayley isomorphic} if there is $\a\in \Aut(G)$ such that $S^{\a}=T$. Cayley isomorphic Cayley (di)graphs are isomorphic, but the converse is not true.
A subset $S$ of $G$ with $1\not\in S$ is said to be a {\em CI-subset} if $\Cay(G,S)\cong \Cay(G,T)$, for some $T\subseteq G$ with $1\not\in T$, implies that they are Cayley isomorphic, and in  this case, $\Cay(G,S)$ is called a {\em CI-digraph}, or a {\em CI-graph} when $S=S^{-1}$. A group $G$ is called a {\em DCI-group} if all Cayley digraphs of $G$ are CI-digraphs, and a {\em CI-group} if all Cayley graphs of $G$ are CI-graphs.

DCI-groups and CI-groups have been widely investigated over last fifty years, and they have been reduced to some special restricted groups in~\cite{DE1,Li2}. However, it is still very  difficult to determine whether a particular group is a DCI-group or a CI-group; see \cite{AN,D2,Feng3,Kov,Mu,MS,N,So,Sp,Sp1} for example. In fact, it is even an open problem to classify dihedral DCI-groups or CI-groups.

Motivated by the \'Ad\'am conjecture (see \cite{Ad}): every finite cyclic group is a CI-group, due to contributions of many researchers like Elspas and Turner~\cite{Elspas}, Djokovi\'c~\cite{Djokovic}, Turner~\cite{Turner}, Babai~\cite{Babai}, Alspach and Parsons~\cite{Alspach}, Godsil~\cite{Godsil} and P\'alfy~\cite{Palfy}, the cyclic DCI-groups and CI-groups were finally classified by Muzychuk~\cite{Muzychuk,M.Muzychuk}, which say that a cyclic group of order $n$ is a DCI-group if and only if $n=mk$ where $m=1,2,4$ and $k$ is odd-square-free, and is a CI-group if and only if either $n=8,9,18$, or $n=mk$ where $m=1,2,4$ and $k$ is odd-square-free. The generalised dihedral group
$\mathrm{Dih}(A)$ over an abelian group $A$ is the group $\langle A,b\ | \ b^2=1, a^b = a^{-1}, \forall a \in A\rangle$, and in particular, if $A$ is a cyclic group of order $n$, $\mathrm{Dih}(A)$ is the dihedral group $\mathrm{D}_{2n}$ of order $2n$. There are also some partial results on dihedral DCI-groups or CI-groups. Babai~\cite{Babai} proved that the dihedral group $\mathrm{D}_{2p}$ for a prime $p$ is a CI-group. Conder and Li~\cite{conder} proved that $\mathrm{D}_{18}$ is a CI-group, and in~\cite{DE1}, it was further proved that $\mathrm{D}_{6p}$~($p$ a prime) is a DCI-group if and only if $p \geq 5$, and is a CI-group if and only if $p \geq 3$. Recently, Dobson et al.~\cite{DE2} proved that if $R$ is a generalised dihedral CI-group, then for every odd prime $p$ the Sylow $p$-subgroup of $R$ has order $p$, or 9, which reduces dihedral DCI-groups to $\mathrm{D}_{2n}$ with $n=2^mk$ and $k$ odd-square-free and dihedral CI-groups to $\mathrm{D}_{18}$ or $\mathrm{D}_{2n}$ with $n=2^mk$ and $k$ odd-square-free.

Let $\Cay(G,S)$ be a Cayley (di)graph of $G$ with respect to $S$. For a given $g\in G$, the right multiplication $R(g): x\mapsto xg$, $\forall x\in G$, is an automorphism of $\Cay(G,S)$, and $R(G):=\{R(g)\ | \ g\in G\}$ is a regular group of automorphisms of $\Cay(G,S)$, called the right regular representation of $G$. A Cayley (di)graph $\Cay(G,S)$ is said to be normal if $R(G)$ is a normal subgroup of $\Aut(\Cay(G,S))$. Normality of Cayley (di)graphs is very important because the automorphism groups of normal Cayley  (di)graph are actually known; see Godsil~\cite{Godsil1} or Proposition~\ref{N_AUT}. Furthermore, the study of normality of Cayley (di)graphs is currently a hot topic in algebraic graph theory, and we refer to~\cite{FLuWX,Feng2,Zhou,FMW,DFZ,DF,PYL,YFZC} for example.

A group $G$ is called a {\em NDCI-group} or a {\em NCI-group} if all normal Cayley digraphs or graphs of $G$ are CI-digraphs or CI-graphs, respectively. Obviously, a DCI-group is a NDCI-group and a CI-group is a NCI-group. Li~\cite{C.H.Li} constructed some normal Cayley digraphs of cyclic groups of $2$-power order which are not CI-graphs, and proposed the following problem: Characterize normal Cayley digraphs which are not CI-graphs. Similar to DCI-groups and CI-groups, a natural problem is to classify finite NDCI-groups and NCI-groups. Recently, Xie et al.~\cite{XFRL} classified cyclic NDCI-groups and NCI-groups, and in this paper, we classify dihedral NDCI-groups and NCI-groups.

\begin{thm}\rm\label{mainth} Let $n\geq 2$ be a positive integer and let $\mathrm{D}_{2n}$ be the dihedral group of order $2n$. Then the following facts are equivalent:
\begin{enumerate}
  \item $\D2n$ is a NDCI-group;
  \item $\D2n$ is a NCI group;
  \item Either $n=2,4$ or $n$ is odd.
\end{enumerate}
\end{thm}

Classification of NDCI-groups and NCI-groups can be helpful for classification of DCI-groups and CI-groups. In fact, by Theorem~\ref{mainth}, together with others, we have the following corollary.

\begin{cor}\rm\label{cor-DCI}
If a dihedral group $\mathrm{D}_{2n}$ of order $2n$ is a $\rm{DCI}$-group then $n=2$ or $n$ is odd-square-free, and if $\mathrm{D}_{2n}$ is a $\rm{CI}$-group then $n=2,9$ or $n$ is odd-square-free.
\end{cor}

We believe that the converse of Corollary~\ref{cor-DCI} is true.

\begin{conj}\rm\label{conjectrue}
A dihedral group $\mathrm{D}_{2n}$ of order $2n$ is a $\rm{DCI}$-group if and only if $n=2$ or $n$ is odd-square-free, and  $\mathrm{D}_{2n}$ is a $\rm{CI}$-group if and only if $n=2,9$ or $n$ is odd-square-free.
\end{conj}

Note that Holt and Royle~\cite[6.2]{DHGR} claimed that $\mathrm{D}_8$ is a CI-group, and this is not true by Corollary~\ref{lem-S-C}, where a Cayley graph of $\mathrm{D}_8$ is constructed and its non-CI-property is also checked by using {\sc Magma}~\cite{magma}.

All the notations and terminologies used in this paper are standard, and for group and graph concepts not defined here, we refer to \cite{Biggs,Dixon,Isaacs,Rotman}.
\section{Preliminaries}

In this section, we give some basic concepts and facts that will be needed later. A group $G$ is called a $p$-group if $|G|$ is a $p$-power for some prime $p$. The following is a basic property for $p$-group; see \cite[Theorem 1.19]{Isaacs}.

\begin{pro}\rm\label{p-group}
Let $G$ be a finite $p$-group and let $N\not=1$ be a normal subgroup of $G$. Let $Z(G)$ be the center of $G$. Then $N\cap Z(G)\not=1$.
\end{pro}

Let $G$ be a group. For $x,y\in G$, we use $[x,y]$ as an abbreviation for the
{\em commutator} $x^{-1}y^{-1}xy$ of $x$ and $y$. The following  proposition is fundamental for   commutators and its proof is straightforward; see \cite[Chapter 4]{Isaacs}.

\begin{pro}\rm\label{commutator}
Let $n$ be a positive integer and let $G$ be a group. Then, for all $x, y, z \in G$, we have:
\begin{itemize}
\item [(1)] $[xy, z]=[x, z]^y[y, z]$ and $[z, xy]=[z, y][z, x]^y$;
\item [(2)] If $[x,y]\in Z(G)$, then $(xy)^n=x^ny^n[y,x]^{\binom{n}{2}}$.
\end{itemize}
\end{pro}

A group $G$ is called {\it metacyclic} if $G$ has a normal cyclic group $N$ such that $G/N$ is cyclic, and the metacyclic group $G$ is called {\it split} if $N$ has a complement in $G$, that is $G$ has a subgroup $H$ such that $G=NH$ and $N\cap H=1$; in this case, $G$ is a semiproduct of $N$ by $H$, and we write $G=N\rtimes H$. For a positive integer $n$, denote by $\ZZ_n$ the additive group of integers modulo $n$, and by $\ZZ_n^*$ the multiplicative
group of $\ZZ_n$ consisting of integers coprime to $n$.

\begin{pro}  \rm\label{metacyclic} Let $p$ be an odd prime and let $G$ be a split metacylic $p$-group such that $N=\langle a\rangle\unlhd G$ and $G=N\rtimes \langle b\rangle$ with $o(a)=p^s\geq p$ and $o(b)=p^t\geq p$. Then every element of order $p$ in $G$ lies in the subgroup $\langle a^{p^{s-1}}\rangle\times\langle b^{p^{t-1}}\rangle$.
\end{pro}

\pf Note that $\langle a^{p^{s-1}}\rangle$ and $\langle b^{p^{t-1}}\rangle$ are the unique subgroups of order $p$ in $N$ and $\langle b\rangle$, respectively. Then $\langle a^{p^{s-1}}\rangle$ is characteristic in $N$, and since $N\unlhd G$, we derive $\langle a^{p^{s-1}}\rangle\unlhd G$. By Propostion~\ref{p-group}, $Z(G)\cap \langle a^{p^{s-1}}\rangle\not=1$, implying $\langle a^{p^{s-1}}\rangle\leq Z(G)$. Since $N\unlhd G$ and $G/N\cong\langle b\rangle$, every element of order $p$ in $G$ lies in the subgroup $N\rtimes\langle b^{p^{t-1}}\rangle$.

Let $c=b^{p^{t-1}}$. Since $N=\langle a\rangle\unlhd G$ and $o(a)=p^s$, we may assume $a^c=a^i$ for some $i\in \ZZ_{p^s}^*$. Then $a=a^{c^p}=a^{i^p}$ and hence $i^p=1$ in $\ZZ_{p^s}^*$, that is, $i=1$ or $i$ is an element of order $p$ in $\ZZ_{p^s}^*$. Since $p$ is odd, $\ZZ_{p^s}^*\cong\ZZ_{p^{s-1}(p-1)}$ either has a trivial Sylow $p$-subgroup ($s=1$), or has a unique subgroup of order $p$ ($s\geq 2$). It follows that  $i=1$ when $s=1$ and $i=kp^{s-1}+1$ for some $0\leq k\leq p-1$ when $s\geq 2$. Thus, $[a,c]=a^{-1}a^c=1$ or $a^{kp^{s-1}}$, of which both belong to $Z(G)$.  By Proposition~\ref{commutator}, for all $j\in\ZZ_{p^s}$ we have $[a^j,c]=[a,c]^j\in Z(G)$ and
$a^jc$ has order $p$ if and only if $a^j$ has order $p$, which implies that every element of order $p$ in $G$ lies in the subgroup $\langle a^{p^{s-1}}\rangle\times\langle b^{p^{t-1}}\rangle$.\qed

Let $\Cay(G,S)$ be a Cayley digraph of a group $G$ with respect to $S$. Then $\Aut(G,S)=\{ \alpha\in \Aut(G)~|~S^{\alpha}=S\}$ is a subgroup of $\Aut(\Cay(G,S))_1$, the stabilizer of $1$ in $\Aut(\Cay(G,S))$. By Godsil~\cite{Godsil1}, the normalizer of $R(G)$ in $\Aut(\Cay(G,S))$ is the semiproduct $R(G)\rtimes \Aut(G,S)$, where $R(g)^\a=R(g^\a)$ for all $g\in G$ and $\a\in\Aut(G,S)$, and by \cite[Propositions 1.3 and 1.5]{Xu1}, we have the following.

\begin{pro}\rm\label{N_AUT} Let $\Cay(G,S)$ be a Cayley digraph of a group $G$ with respect to $S$ and let $A=\Aut(\Cay(G,S))$. Then $N_A(R(G))=R(G)\rtimes\Aut(G,S)$ and
$\Cay(G,S)$ is normal if and only if $A_1=\Aut(G,S)$.
\end{pro}

Babai~\cite{Babai} gave a well-known criterion for a Cayley digraph to be a CI-digraph, that is, a Cayley digraph $\Cay(G,S)$ is a CI-digraph if and only if every regular group of automorphisms of $\Cay(G,S)$ isomorphic to $G$, is conjugate to $R(G)$ in $\Aut(\Cay(G,S))$. Based on this, we have the following proposition (also see \cite[Corollary 6.9]{C.H.Li}).

\begin{pro}\rm\label{CI-graph-prop}
Let $\Cay(G,S)$ be a normal Cayley digraph of a group $G$ with respect to $S$. Then $\Cay(G,S)$ is a $\rm{CI}$-digraph if and only if $\Aut(\Cay(G,S))$ has a unique regular subgroup isomorphic to $G$, that is, $R(G)$.
\end{pro}

The following proposition was  given in \cite[Proposition~2.6]{XFRL}, which is about a family of non-normal Cayley digraphs, called generalised wreath product of Cayley digraphs.

\begin{pro} \rm\label{General-normal}  Let $G$ be a finite group, and let $1\not\in S\subseteq G$, $1<H\leq K<G$ and $H\unlhd G$. Assume that $S\backslash K$ is a union of some cosets of $H$ in $G$ and that there exist $x\not\in K$ and $y\in H$ such that $y^x\not=y^{-1}$. Then the Cayley digraph  $\Cay(G,S)$ is non-normal.
\end{pro}


\section{Automorphisms and Holomorphs of Dihedral groups}

In this section we collect some details about automorphism groups and holomorphs of dihedral groups, which are crucial for the proof of Theorem~\ref{mainth}.

Let $G$ be a finite group and let $g\in G$. Denote by $o(g)$ the order of $g$ in $G$. Let $p$ be a prime and $\pi$ a set of primes. Denote by $G_p$ a Sylow $p$-subgroup of $G$. An element $g$ of $G$ is called a {\em $\pi$-element} if all prime factors of $o(g)$ belong to $\pi$, and a {\em $p'$-element} if $o(g)$ has no factor $p$. If $G$ is soluble, denote by $G_\pi$ a Hall $\pi$-subgroup of $G$, and by $G_{p'}$ a Hall $p'$-subgroup of $G$.

Let $n$ be a positive integer. We first make the following convention:
\begin{equation}\label{eq1}
n={\prod}_{i=1}^m{p_i^{r_i}}\geq 2, \ C_n=C_{p_1^{r_1}}\times\ldots\times C_{p_m^{r_m}}=\langle a_1\rangle\times\ldots\times\langle a_m\rangle,\ a=a_1\ldots a_m,
\end{equation}
where $p_1, \ldots, p_m$ are all distinct prime factors of $n$, $C_{p_i^{r_i}}=\langle a_i\rangle\cong \ZZ_{p_i^{r_i}}$ for each $1\leq i\leq m$, and $C_n=\langle a\rangle\cong \ZZ_n$. By \cite[Theorem 7.3]{Rotman}, we have
\begin{equation}\label{eq2}
\Aut(C_n)=\Aut(C_{p_1^{r_1}})\times \Aut(C_{p_2^{r_2}})\times\ldots\times\Aut(C_{p_m^{r_m}}),
\end{equation}
where $\Aut(C_{p_i^{r_i}})$ is viewed as the subgroup of $\Aut(C_n)$ by identifying $\a_i\in \Aut(C_{p_i^{r_i}})$ as the automorphism of $C_n$ induced by $a_i\mapsto a_i^{\a_i}$ and $a_j\mapsto a_j$ for $j\not=i$. Further we set
\begin{equation}\label{eq3}
\D2n=\langle a,b\ |\ a^n=b^2=1, bab=a^{-1}\rangle=C_n\rtimes\langle b\rangle.
\end{equation}
For $a^i\in \langle a\rangle$, denote by $\theta_{a^i}$ the automorphism of $\D2n$ induced by
\begin{equation}\label{eq4}
\theta_{a^i}: a\mapsto a, b\mapsto ba^i.
\end{equation}
Then $o(\theta_{a^i})=o(a^i)$ and $\langle \theta_{a}\rangle\cong \ZZ_n$. Write $\overline{a}_i=a_1\ldots a_{i-1}a_{i+1}\ldots a_m$. Then $o(\overline{a}_i)=n/p_i^{r_i}$, and hence
 $\langle \theta_{a}\rangle=\langle \theta_{a_i}\rangle\times\langle \theta_{\overline{a}_i}\rangle$. In fact, $\langle \theta_{a_i}\rangle$ and $\langle \theta_{\overline{a}_i}\rangle$ are the unique Sylow $p_i$-subgroup and the unique Hall $p_i'$-subgroup of $\langle \theta_{a}\rangle$, respectively. Clearly,
 \begin{equation}\label{eq5}
\theta_{a^i}\theta_{a^j}=\theta_{a^{i+j}} \mbox{ and } \theta_{a^i}^k=\theta_{a^{ik}} \mbox{ for all  $i,j,k\in \ZZ_n$}.
\end{equation}
In what follows, we also view $\Aut(C_n)$ as the subgroup of $\Aut(\D2n)$ by identifying $\b\in \Aut(C_n)$ as the automorphism of $\Aut(\D2n)$ induced by $a\mapsto a^\b$ and $b\mapsto b$. Then for each $1\leq i\leq m$, we have $\Aut(C_{p_i^{r_i}})\leq \Aut(\D2n)$.

\begin{lem}  \rm\label{AutD2n} Let $n\geq 3$. Then $\Aut(\D2n)$ has the following properties.
\begin{itemize}
\item [(1)] $\Aut(\D2n)=\langle \theta_{a}\rangle\rtimes \Aut(C_n)$, and $\langle \theta_{a}\rangle$ is the kernel of the natural action of $\Aut(\D2n)$ on $C_n$. Furthermore,  $\theta_{x}^\b=\theta_{x^\b}$ for all $x\in\langle a\rangle$ and $\b\in \Aut(C_n)$;

\item [(2)] For every $1\leq i\leq m$, $\langle \theta_{a}\rangle\rtimes \Aut(C_{p_i^{r_i}})=(\langle \theta_{a_i}\rangle\rtimes \Aut(C_{p_i^{r_i}})) \times \langle \theta_{\overline{a}_i}\rangle$, where $\overline{a}_i=a_1\ldots a_{i-1}a_{i+1}\ldots a_m$;

\item [(3)] Assume $p_1>p_2>\ldots>p_m$. For every $1\leq i\leq m$, set $\pi_i=\{p_1,\ldots,p_i\}$. Then $\Aut(\D2n)$ has a normal Hall $\pi_i$-subgroup and $\Aut(\D2n)_{\pi_i}=\langle \theta_{a_1\ldots a_i}\rangle\rtimes \Aut(C_n)_{\pi_i}$, where $\Aut(C_n)_{\pi_i}=\Aut(C_{p_1^{r_1}})_{\pi_i}\times\ldots\times \Aut(C_{p_i^{r_i}})_{\pi_i}$;

    \item [(4)] Let $p_i$ be odd and let $\g\in\langle \theta_{a}\rangle\rtimes \Aut(C_{p_i^{r_i}})$ with $o(\g)=p_i$. Then in $\D2n$, $\g$ fixes each coset of $H=\langle a^{n/p_i}\rangle$, and if $\g\not\in \langle \theta_{a}\rangle$ then $r_i\geq 2$ and $\g$ has fixed-point-set $M\cup ba^rM$ for some $1\leq r\leq p_i-1$, where $M=\langle a_1,\ldots,a_{i-1},a_i^{p_i},a_{i+1},\ldots,a_m\rangle$.
\end{itemize}

\end{lem}

\pf Since $n\geq 3$, $C_n$ is the unique cyclic subgroup of order $n$ in $\D2n$, and hence $C_n$ is characteristic in $\D2n$, implying that $\Aut(\D2n)$ fixes $C_n$ setwise. Let $K$ be the kernel of $\Aut(\D2n)$ acting on $C_n$. Then for all $\g\in K$ we have $a^\g=a$, and $\Aut(\D2n)/K$ acts faithfully on $C_n$. Thus, $|\Aut(\D2n)/K|\leq |\Aut(C_n)|$. Since $\Aut(C_n)$ is viewed as a subgroup of $\Aut(\D2n)$ fixing $b$, we have $K\cap \Aut(C_n)=1$ and $\Aut(\D2n)=K\rtimes\Aut(C_n)$.

For every $ba^i$ with $0\leq i< n$, $\D2n$ has an automorphism induced by $a\mapsto a$ and $b\mapsto ba^i$, and only the identity in $\Aut(\D2n)$ fixes $a$ and $ba^i$. Therefore, $|K|=n$. Since $\theta_{a}\in K$ and $\langle \theta_{a}\rangle\cong \ZZ_n$, we have $K=\langle\theta_a\rangle$, and then $\Aut(\D2n)=\langle \theta_{a}\rangle\rtimes \Aut(C_n)$.

For $\b\in\Aut(C_n)$, we have $b^\b=b$, and for $x\in\langle a\rangle$, $\theta_{x}$ fixes $\langle a\rangle$ pointwise. Since $\b$ fixes $\langle a\rangle$ setwise, $$a^{{\theta_{x}}^\b}=(a^{\b^{-1}})^{\theta_{x}\b}=a^{\b^{-1}\b}=a=a^{\theta_{x^\b}},$$ and $$b^{{\theta_{x}}^\b}=(b^{\b^{-1}})^{\theta_{x}\b}=b^{\theta_{x}\b}=(bx)^\b=bx^\b=b^{\theta_{x^\b}}.$$ Since $\D2n=\langle a,b\rangle$, we obtain $\theta_{x}^\b=\theta_{x^\b}$. This completes the proof of part~(1).

Recall that $\langle \theta_{a}\rangle=\langle \theta_{a_i}\rangle\times\langle \theta_{\overline{a}_i}\rangle$. Let $\b\in \Aut(C_{p_i^{r_i}})$. Since $\overline{a}_i=a_1\ldots a_{i-1}a_{i+1}\ldots a_m$, we have $\overline{a}_i^\b=\overline{a}_i$, and $\theta_{\overline{a}_i}^\b=\theta_{(\overline{a}_i)^\b}=\theta_{\overline{a}_i}$, that is, $\theta_{\overline{a}_i}\b=\b\theta_{\overline{a}_i}$. Thus, $$\langle \theta_{a}\rangle\rtimes \Aut(C_{p_i^{r_i}})=(\langle \theta_{a_i}\rangle\rtimes \Aut(C_{p_i^{r_i}})) \times \langle \theta_{\overline{a}_i}\rangle.$$  This completes the proof of part~(2).

Let $1\leq i\leq m$. Note that $\Aut(C_{p_i^{r_i}})\cong \ZZ_{p_i^{r_i-1}(p_i-1)}$ if $p_i$ is odd, and $\Aut(C_{p_i^{r_i}})\cong \ZZ_{p_i}\times \ZZ_{p_i^{r_i-2}}$ if $p_i=2$ and $r_i\geq 2$. Since $p_1>p_2>\ldots>p_m$, we have $\Aut(C_{p_k^{r_k}})_{\pi_i}=1$ for $i<k\leq m$, and hence $\Aut(C_n)_{\pi_i}=\Aut(C_{p_1^{r_1}})_{\pi_i}\times\ldots\times \Aut(C_{p_i^{r_i}})_{\pi_i}$.

It is easy to see that $\theta_{a_1\ldots a_i}=\theta_{a_1}\ldots\theta_{a_i}$ has order $p_1^{r_1}\ldots p_i^{r_i}$. Then $\langle \theta_{a_1\ldots a_i}\rangle$ is a Hall $\pi_i$-subgroup of $\langle\theta_a\rangle$, and hence $\langle \theta_{a_1\ldots a_i}\rangle$ is characteristic in $\langle\theta_a\rangle$. Since $\langle\theta_a\rangle\unlhd \Aut(\D2n)$, we have  $\langle \theta_{a_1\ldots a_i}\rangle\unlhd \Aut(\D2n)$, and $\langle \theta_{a_1\ldots a_i}\rangle\rtimes \Aut(C_n)_{\pi_i}$ is a Hall $\pi_i$-subgroup of $\Aut(\D2n)$.

Recall that $\Aut(\D2n)=\langle \theta_{a}\rangle\rtimes \Aut(C_n)$. To prove $\langle \theta_{a_1\ldots a_i}\rangle\rtimes \Aut(C_n)_{\pi_i}\unlhd \Aut(\D2n)$, it suffices to show that $\Aut(C_n)_{\pi_i}^{\theta_{a}}\leq \langle \theta_{a_1\ldots a_i}\rangle\rtimes \Aut(C_n)_{\pi_i}$, or alternatively, $\Aut(C_{p_j^{r_j}})_{\pi_i}^{\theta_{a}}\leq \langle \theta_{a_1\ldots a_i}\rangle\rtimes \Aut(C_n)_{\pi_i}$ for every $1\leq j\leq i$. Now take $\a\in \Aut(C_{p_j^{r_j}})_{\pi_i}$. Then $a_k^\a=a_k$ for $1\leq k\leq m$ with $k\not=j$, and $a_j^\a\in \langle a_j\rangle$. It follows that  $$\a^{\theta_a}=\theta_a^{-1}\a\theta_a=\theta_a^{-1}\theta_a^{\a^{-1}}\a=
\theta_{a^{-1}}\theta_{a^{\a^{-1}}}\a=\theta_{a^{-1}a^{\a^{-1}}}\a=\theta_{a_j^{-1}a_j^{\a^{-1}}}\a\in
\langle \theta_{a_j}\rangle\Aut(C_{p_j^{r_j}})_{\pi_i},$$ where $$a^{-1}a^{\a^{-1}}=a_1^{-1}a_2^{-1}\ldots a_m^{-1}a_1\ldots a_{j-1}a_j^{\a^{-1}}a_{j+1}\ldots a_m=a_j^{-1}a_j^{\a^{-1}}\in \langle a_j\rangle.$$ Since $1\leq j\leq i$, we have $$\langle \theta_{a_j}\rangle\Aut(C_{p_j^{r_j}})_{\pi_i}\leq \langle \theta_{a_1\ldots a_i}\rangle\rtimes \Aut(C_n)_{\pi_i},$$ and hence $\Aut(C_{p_j^{r_j}})_{\pi_i}^{\theta_{a}}\leq \langle \theta_{a_1\ldots a_i}\rangle\rtimes \Aut(C_n)_{\pi_i}$, as required. Thus, $\Aut(\D2n)$ has a normal Hall $\pi_i$-subgroup, that is, $\Aut(\D2n)_{\pi_i}=\langle \theta_{a_1\ldots a_i}\rangle\rtimes \Aut(C_n)_{\pi_i}$. This complete the proof of part~(3).

Let $p_i$ be odd, and let $\a\in\langle \theta_{a}\rangle$ such that $a^\a=a$ and $b\mapsto ba_i^{p_i^{r_i-1}}$, that is, $\a=\theta_{a_i^{p_i^{r_i-1}}}$. Then $\langle \a\rangle$ is the unique subgroup of order $p_i$ in $\langle \theta_a\rangle$. Since both $\langle a^{n/p_i}\rangle$ and $\langle a_i^{p_i^{r_i-1}}\rangle$ are the  unique subgroup of order $p_i$ in $\D2n$, we have $H=\langle a^{n/p_i}\rangle=\langle a_i^{p_i^{r_i-1}}\rangle$. Clearly, $\a$ fixes each coset of $H$ in $\D2n$ and has no fixed point in $b\langle a\rangle$. Since $\D2n=\langle a\rangle\cup b\langle a\rangle$, $\a$ has fixed-point-set $\langle a\rangle$ in $\D2n$.

Since $p_i$ is odd, $\Aut(C_{p_i^{r_i}})\cong\ZZ_{p_i^{r_i-1}(p_i-1)}$, and hence $\Aut(C_{p_i^{r_i}})$ has an element of order $p_i$ if and only if $r_i\geq 2$. Note that $\g\in\langle \theta_{a}\rangle\rtimes \Aut(C_{p_i^{r_i}})$ with $o(\g)=p_i$.
If $r_i=1$ then $\Aut(C_{p_i^{r_i}})$ has no element of order $p_i$, and in this case, we may let  $\g=\a\in\langle \theta_{a}\rangle\rtimes \Aut(C_{p_i^{r_i}})$. Then $\g$ fixes each coset of $H$ in $\D2n$.

Assume $r_i\geq 2$. Let $\b\in \Aut(C_{p_i^{r_i}})$ such that $a_i^\b=a_ia_i^{p_i^{r_i-1}}$, $b^\b=b$ and $a_j^\b=a_j$ for $j\not=i$. Then $a^\b=aa_i^{p_i^{r_i-1}}$, and $\langle \b\rangle$ is the unique subgroup of order $p_i$ in $\Aut(C_{p_i^{r_i}})$. Recall that $H=\langle a_i^{p_i^{r_i-1}}\rangle$. Clearly, $\b$ fixes each coset of $H$ in $\D2n$, and fixes $$M=\langle a_1,\ldots,a_{i-1},a_i^{p_i},a_{i+1},\ldots,a_m\rangle \mbox{ and }\ bM, $$
pointwise. For every $1\leq t\leq p_i-1$, it is easy to see that $(a_i^t)^\b=a_i^ta_i^{tp_i^{r_i-1}}\not=a_i^t$ and $(ba_i^t)^\b\not=ba_i^t$, so $\b$ has no fixed point in the coset $a_i^tM$ or $ba_i^tM$. Since  $\D2n=\cup_{k=0}^{p_i-1} a^kM\cup ba^kM$, $M\cup bM$ is the fixed-point-set of $\b$ in $\D2n$.

By part~(2), $\langle \theta_{a}\rangle\rtimes \Aut(C_{p_i^{r_i}})$ has a normal Sylow $p_i$-subgroup that is a split metacyclic $p_i$-group, and by Proposition~\ref{metacyclic}, $\g\in\langle \a\rangle\times\langle \b\rangle$. Since both $\a$ and $\b$ fix each coset of $H$, $\g$ fixes each coset of $H$ in $\D2n$.

Recall that $\g\in\langle \theta_{a}\rangle\rtimes \Aut(C_{p_i^{r_i}})$ with $o(\g)=p_i$. To finish the proof of (4), let $\g\not\in \langle \theta_{a}\rangle$. Then $\Aut(C_{p_i^{r_i}})$ has an element of order $p_i$, which forces $r_i\geq 2$. Since $\g\in\langle \a\rangle\times\langle \b\rangle$, we have $\g=\a^k\b^\ell$ for some $k,\ell\in \ZZ_{p_i}$ with $\ell\not=0$. Then $\b^\ell\not=1$ and $\langle \b^\ell\rangle=\langle \b\rangle$ as $o(\b)=p_i$. Since $x^\g=x^{\a^k\b^\ell}=x^{\b^\ell}$ for every $x\in \langle a\rangle$, $\g$ and $\b^\ell$ has the same fixed-point-set in $\langle a\rangle$. Thus, $\g$ has fixed-point-set $M$ in $\langle a\rangle$, because $\b$ has fixed-point-set $M$ in $\langle a\rangle$.

Now consider the fixed-point-set of $\g$ in $b\langle a\rangle$. Since $\ell\not=0$, $\ell$ has an inverse in the field $\ZZ_{p_i}$ of order $p_i$, say $\ell^{-1}$. Take $r=-k\ell^{-1}$ in $\ZZ_{p_i}$. Note that $r_i\geq 2$ and $a=a_1\ldots a_m$. Then $$(ba^r)^\g=(ba^r)^{\a^k\b^\ell}=(ba_i^{kp_i^{r_i-1}}a^r)^{\b^\ell}=b(a_i^{kp_i^{r_i-1}}a^r)a_i^{r\ell p_i^{r_i-1}}=ba^r,$$ which implies that $\g$ fixes $ba^rM$ pointwise. If $(ba^t)^\g=ba^t$ for some $1\leq t\leq p_i-1$ with $t\not=r$, then $\g$ fixes $a^{t-r}$ and so $a^{t-r}M$ pointwise,  contradicting that $\g$ has fixed-point set $M$ in $\langle a\rangle$. Thus, $\g$ has fixed-point-set $ba^rM$ in $b\langle a\rangle$ and hence $\g$ has fixed-point-set $M\cup ba^rM$ in $\D2n$. This completes the proof of part~(4). \qed

For a finite group $G$, the right regular representation $R(G)$ and the automorphism group $\Aut(G)$ are permutation groups on $G$. Furthermore, $R(G)\cap \Aut(G)=1$, and $R(G)\Aut(G)=R(G)\rtimes\Aut(G)$, where $R(g)^\a=R(g^\a)$ for all $g\in G$ and $\a\in \Aut(G)$. The normalizer of $R(G)$ in the symmetric group $S_G$ on $G$ is called the {\em holomorph} of $G$, denoted by $\Hol(G)$, and by \cite[Lemma 7.16]{Rotman}, $\Hol(G)=R(G)\rtimes\Aut(G)$. Now we have
\begin{equation}\label{eq6}
\Hol(\D2n)=R(\D2n)\rtimes\Aut(\D2n)=R(\D2n)\rtimes (\langle \theta_{a}\rangle\rtimes\Aut(C_n)).
\end{equation}
Note that $\Aut(C_n)=\Aut(C_{p_1^{r_1}})\times \Aut(C_{p_2^{r_2}})\times\ldots\times\Aut(C_{p_m^{r_m}})$ by Eq~(\ref{eq2}).

\begin{lem}  \rm\label{holomorph}
Let $n$ be odd. Using the notations and formulae in Eqs~(\ref{eq6})-(\ref{eq6}), we have the following.
\begin{itemize}
\item [(1)] For all $d\in \D2n$ and $\a\in\Aut(\D2n)$, we have $R(d)^\a=R(d^\a)$, and $\langle R(a)\rangle\langle\theta_a\rangle=\langle R(a)\rangle\times\langle\theta_a\rangle\unlhd \Hol(\D2n)$;
\item [(2)] Assume $p_1>p_2>\ldots>p_m$. For each $1\leq i\leq m$, set $\pi_i=\{p_1,\ldots,p_i\}$. Then  $\Hol(\D2n)_{\pi_i}\leq(\langle R(a)\rangle\times\langle\theta_a\rangle)\rtimes \Aut(C_n)_{\pi_i}\unlhd \Hol(\D2n)$, where $\Aut(C_n)_{\pi_i}=\Aut(C_{p_1^{r_1}})_{\pi_i}\times\ldots\times \Aut(C_{p_i^{r_i}})_{\pi_i}$. Furthermore, $\langle R(b)\rangle\Aut(C_n)_2=\langle R(b)\rangle\times\Aut(C_n)_2$ is a Sylow $2$-subgroup of $\Hol(\D2n)$, where $\Aut(C_n)_2=\Aut(C_{p_1^{r_1}})_2\times \Aut(C_{p_2^{r_2}})_2\times\ldots\times \Aut(C_{p_m^{r_m}})_2$;
\item [(3)] Let $1\leq i,j\leq m$ with $i\not=j$. Then $\Aut(C_{p_i^{r_i}})$, under conjugacy, fixes each $p_j$-element in $\langle R(a)\rangle\times\langle\theta_a\rangle$, and
if $\a\in \Aut(C_{p_i^{r_i}})$ fixes an element of order $p_i^{r_i}$ in $\langle R(a)\rangle\times \langle \theta_a\rangle$, then $\a=1$.
\end{itemize}
\end{lem}

\pf It is known that $R(d)^\a=R(d^\a)$ for all $d\in \D2n$ and $\a\in\Aut(\D2n)$. Then $R(a)^{\theta_a}=R(a^{\theta_a})=R(a)$, that is, $R(a)$ commutes with $\theta_a$. Since $\langle R(a)\rangle\cap\langle\theta_a\rangle=1$, we have
$\langle R(a)\rangle\langle\theta_a\rangle=\langle R(a)\rangle\times\langle\theta_a\rangle$.
 Clearly, $\langle R(a)\rangle\times\langle\theta_a\rangle$ is a Hall $2'$-subgroup of $R(\D2n)\rtimes \langle \theta_a\rangle$, and hence $\langle R(a)\rangle\times\langle\theta_a\rangle$ is characteristic in $R(\D2n)\rtimes \langle \theta_a\rangle$. It follows that $\langle R(a)\rangle\times\langle\theta_a\rangle\unlhd \Hol(\D2n)$, as $R(\D2n)\rtimes \langle \theta_a\rangle\unlhd \Hol(\D2n)$ by Eq~(\ref{eq6}). This completes the proof of part~(1).

To prove part~(2), assume $p_1>p_2>\ldots>p_m$, and  for each $1\leq i\leq m$, let $\pi_i=\{p_1,\ldots,p_i\}$. By part~(1), $\langle R(a)\rangle\times\langle\theta_a\rangle\unlhd \Hol(\D2n)$, and by Eq~(\ref{eq6}),
$$\Hol(\D2n)=R(\D2n)\rtimes (\langle \theta_{a}\rangle\rtimes\Aut(C_n))=(R(\D2n)\rtimes \langle \theta_{a}\rangle)\rtimes\Aut(C_n).$$
Write $A=\langle R(a)\rangle\times\langle\theta_a\rangle$.
Then $A\cap \Aut(C_n)=1$ and $(R(\D2n)\rtimes \langle \theta_{a}\rangle)/A$ is a normal subgroup of order $2$ in $\Hol(\D2n)/A$, and hence lies in the center of $\Hol(\D2n)/A$. Thus,
$$\Hol(\D2n)/A=(R(\D2n)\rtimes \langle \theta_{a}\rangle)/A\times \Aut(C_n)A/A,$$ where
$\Aut(C_n)A/A\cong \Aut(C_n)$. Since $\Aut(C_n)$ is abelian, $\Hol(\D2n)/A$ is abelian, and therefore, $(\langle R(a)\rangle\times\langle\theta_a\rangle)\rtimes \Aut(C_n)_{\pi_i}=A\rtimes \Aut(C_n)_{\pi_i}\unlhd \Hol(\D2n)$.

Since $n$ is odd, all $p_i$ are odd and hence $2\not\in \pi_i$ for each $1\leq i\leq m$. Clearly, $$\Hol(\D2n)/((R(\D2n)\rtimes\langle\theta_a\rangle)\rtimes \Aut(C_n)_{\pi_i})\cong \Aut(C_n)/\Aut(C_n)_{\pi_i}.$$ Then for every Hall $\pi_i$-subgroup $\Hol(\D2n)_{\pi_i}$,
$\Hol(\D2n)_{\pi_i}\leq (R(\D2n)\rtimes\langle\theta_a\rangle)\rtimes \Aut(C_n)_{\pi_i}$, and it follows that $\Hol(\D2n)_{\pi_i}\leq \langle R(a)\rangle\times\langle\theta_a\rangle)\rtimes \Aut(C_n)_{\pi_i}$, as $|(R(\D2n)\rtimes\langle\theta_a\rangle):(\langle R(a)\rangle\times\langle\theta_a\rangle)|=2$.
By  Lemma~\ref{AutD2n}~(3), $\Aut(C_n)_{\pi_i}=\Aut(C_{p_1^{r_1}})_{\pi_i}\times\ldots\times \Aut(C_{p_i^{r_i}})_{\pi_i}$.

Since $\Aut(C_n)$ fixes $b$, we have $\langle R(b)\rangle\Aut(C_n)=\langle R(b)\rangle\times \Aut(C_n)$, and  therefore, $\langle R(b)\rangle\Aut(C_n)_2=\langle R(b)\rangle\times \Aut(C_n)_2$. By Eq~(\ref{eq6}), $|\Hol(\D2n)_2|=2|\Aut(C_n)_2|$, so $\langle R(b)\rangle\times \Aut(C_n)_2$ is a Sylow $2$-subgroup of $\Hol(\D2n)$. Clearly, $\Aut(C_n)_2=\Aut(C_{p_1^{r_1}})_2\times \Aut(C_{p_2^{r_2}})_2\times\ldots\times \Aut(C_{p_m^{r_m}})_2$. This completes the proof of part~(2).

To prove part~(3), let $1\leq i\leq m$  with $j\not=i$. Recall that $p_i$ is odd. By part~(1), $\langle R(a)\rangle\times\langle\theta_a\rangle\unlhd \Hol(\D2n)$.

Let $\a\in \Aut(C_{p_i^{r_i}})$ and let $x$ be a $p_j$-element in $\langle R(a)\rangle\times \langle \theta_a\rangle$. Since $\langle R(a_j)\rangle \times \langle \theta_{a_j}\rangle$ is a normal Sylow $p_j$-subgroup of $\langle R(a)\rangle\times \langle \theta_a\rangle$, we may let $x=R(a_j)^s\theta_{a_j}^t$ for some $s,t\in\ZZ_{p_j^{r_j}}$. Since $j\not=i$, $a_j^\a=a_j$ and hence $$x^\a=(R(a_j)^s\theta_{a_j}^t)^\a=R(a_j^\a)^s\theta_{a_j^\a}^t=R(a_j)^s\theta_{a_j}^t=x.$$ Thus, $\Aut(C_{p_r^{r_i}})$ fixes each $p_j$-element in $\langle R(a)\rangle\times\langle\theta_a\rangle$.

Now assume that $\a\in \Aut(C_{p_i^{r_i}})$ fixes an element of order $p_i^{r_i}$ in $\langle R(a)\rangle\times \langle \theta_a\rangle$. Since $\langle R(a_i)\rangle \times \langle \theta_{a_i}\rangle$ is a normal Sylow $p_i$-subgroup of $\langle R(a)\rangle\times \langle \theta_a\rangle$, $\a$ fixes an element of order $p_i^{r_i}$ in $\langle R(a_i)\rangle\times \langle \theta_{a_i}\rangle$, say $R(a_i)^k\theta_{a_i}^\ell$ for some $k,\ell\in \ZZ_{p_i^{r_i}}$. Since $o(R(a_i))=o(\theta_{a_i})=p_i^{r_i}$, we have $(k,p_i)=1$ or $(\ell,p_i)=1$. Furthermore, $$R(a_i^k)\theta_{a_i^\ell}=R(a_i)^k\theta_{a_i}^\ell=(R(a_i)^k\theta_{a_i}^\ell)^{\a}=
R((a_i^k)^{\a})\theta_{(a_i^\ell)^{\a}},$$ and since $\langle R(a)\rangle\cap \langle \theta_{a}\rangle=1$, we have $(a_i^k)^{\a}=a_i^k$ and
$(a_i^\ell)^{\a}=a_i^\ell$, which imply that $(a_i)^{\a}=a_i$ as $(k,p_i)=1$ or $(\ell,p_i)=1$. Thus, $\a=1$, completing the proof of part~(3). \qed

\section{Proof of Theorem~\ref{mainth}}

First we consider non-normal Cayley digraphs of a dihedral group admitting a special automorphism of prime order of the digraph. Note that in this section we still use the notations and formulae given in Eqs~(\ref{eq1})-(\ref{eq6}) without explanation.

\begin{lem}  \rm\label{non-normal}
Let $\Cay(\mathrm{D}_{2n},S)$ be a Cayley digraph of the dihedral group $\mathrm{D}_{2n}$ of order $2n$. Assume that $p_t$ is an odd prime for some $1\leq t\leq m$ and $\Aut(\mathrm{D}_{2n},S)$ contains an element of order $p_t$ in $\langle \theta_a\rangle\rtimes\Aut(\bbZ_{p_t^{r_t}})$ but not in $\langle\theta_a\rangle$. Then $\Cay(\mathrm{D}_{2n},S)$ is non-normal.
\end{lem}

\pf Write $\Gamma=\Cay(\mathrm{D}_{2n},S)$ and $A=\Aut(\Gamma)$. Recall that
$n={\prod}_{i=1}^m{p_i^{r_i}}$, $C_n=C_{p_1^{r_1}}\times\ldots\times C_{p_m^{r_m}}=\langle a_1\rangle\times\ldots\times\langle a_m\rangle=\langle a\rangle$ with $o(a_i)=p_i^{r_i}$ and $a=a_1a_2\ldots a_m$, and
 $\D2n=\langle a,b\ |\ a^n=b^2=1, bab=a^{-1}\rangle=C_n\rtimes\langle b\rangle$. By assumption, we may let $z\in \Aut(\D2n,S)\leq \Aut(\Ga)$, $o(z)=p_t$,
 $z\in \langle \theta_a\rangle\rtimes\Aut(\bbZ_{p_t^{r_t}})$ and $z\not\in \langle \theta_a\rangle$. Set
$$M=\langle a_1,\ldots,a_{t-1},a_t^{p_t},a_{t+1},\ldots,a_m\rangle \ \ \mbox{ and } \ \ H=\langle a_t^{p_t^{r_t-1}}\rangle.$$
Then $|C_n:M|=|H|=p_t$. By Lemma~\ref{AutD2n}~(4), $r_t\geq 2$ and hence $1<H\leq M< C_n<\mathrm{D}_{2n}$. Furthermore, in $\D2n$, $z$ has fixed-point-set $M\cup ba^rM$ for some $1\leq r\leq p_t-1$, and $\langle z\rangle$ is transitive on $cH$ for every $c \in \mathrm{D}_{2n}\backslash (M\cup ba^rM)$. Let $$K=M\cup ba^rM,  \ x=a_t^{p_t^{r_t-1}},  \ \mbox{ and } \ y=a_t.$$

Since $C_n$ is characteristic in $\D2n$, we have $M\unlhd \D2n$, and since $o(ba^r)=2$, we have   $K=M\rtimes\langle ba^r\rangle< \D2n$. Thus, $1<H< K<\mathrm{D}_{2n}$. Clearly, $x\in H$, $y\not\in K$ and $x^y\not=x^{-1}$ as $p_t$ is odd. By Proposition~\ref{General-normal}, to prove the lemma, we are left to show that $S\backslash K$ is a union of some cosets of $H$ in $\D2n$. Let $d\in S$ and $d\not\in K$. Since $z\in \Aut(\D2n,S)\leq \Aut(\Ga)$ fixes $K$ pointwise and is transitive on $dH$, we have $dH\subseteq S$, and hence $S\backslash K$ is a union of some cosets of $H$ in $\D2n$. This completes the proof. \qed

\medskip

Next we construct normal Cayley graphs on dihedral groups which are not CI-graphs.
\begin{lem}\rm\label{lem-S}
Let $n>4$ be even. Then  $\Cay(\mathrm{D}_{2n},\{a,a^{-1},b\})$ is a normal Cayley graph, but not a CI-graph.
\end{lem}

\pf By Eq~(\ref{eq3}), $\D2n=\langle a,b\ |\ a^n=b^2=1, bab=a^{-1}\rangle$. Let $\Sigma=\Cay(\mathrm{D}_{2n},S)$ with $S=\{a,a^{-1},b\}$, and let $A=\Aut(\Sigma)$. Since $S=S^{-1}$ and $\langle S\rangle=\mathrm{D}_{2n}$, $\Sigma$ is a connected graph of order $2n$. Write $n=2^sm$ with $(2,m)=1$. Then $s\geq 1$, and one may depict $\Sigma$ as the following figure.

\begin{figure}[htb]
\begin{center}
\begin{picture}(280,135)(10,20)
\thicklines
\qbezier(10,80)(15,135)(70,140)
\qbezier(10,80)(15,25)(70,20)
\qbezier(50,80)(55,95)(70,100)
\qbezier(50,80)(55,65)(70,60)
\put(10,80){\line(1,0){40}}
\put(70,100){\line(1,0){40}}
\put(70,140){\line(1,0){40}}
\put(110,100){\line(1,0){10}}
\put(110,140){\line(1,0){10}}
\put(70,100){\line(0,1){40}}
\put(110,100){\line(0,1){40}}
\put(70,60){\line(0,-1){40}}
\put(70,60){\line(1,0){40}}
\put(70,20){\line(1,0){40}}
\put(110,20){\line(1,0){10}}
\put(110,60){\line(1,0){10}}
\put(110,60){\line(0,-1){40}}
\qbezier(220,140)(275,135)(280,80)
\qbezier(220,20)(275,25)(280,80)
\qbezier(220,100)(235,95)(240,80)
\qbezier(220,60)(235,65)(240,80)
\put(180,140){\line(-1,0){10}}
\put(180,140){\line(1,0){40}}
\put(180,100){\line(-1,0){10}}
\put(180,100){\line(1,0){40}}
\put(180,100){\line(0,1){40}}
\put(220,100){\line(0,1){40}}
\put(180,60){\line(-1,0){10}}
\put(180,60){\line(1,0){40}}
\put(180,20){\line(-1,0){10}}
\put(180,20){\line(1,0){40}}
\put(180,20){\line(0,1){40}}
\put(220,20){\line(0,1){40}}
\put(240,80){\line(1,0){40}}
\multiput(0,20)(0,40){4}{\multiput(70,0)(40,0){2}{\circle*{5}}}
\multiput(0,20)(0,40){4}{\put(180,0){\circle*{5}}}
\multiput(0,20)(0,40){4}{\put(220,0){\circle*{5}}}
\multiput(0,20)(0,40){4}{\put(132,0){\circle*{3}}}
\multiput(0,20)(0,40){4}{\put(145,0){\circle*{3}}}
\multiput(0,20)(0,40){4}{\put(158,0){\circle*{3}}}
\put(10,80){\circle*{5}}
\put(50,80){\circle*{5}}
\put(280,80){\circle*{5}}
\put(240,80){\circle*{5}}
\put(0,77){${\small 1}$}
\put(43,84){${\small b}$}
\put(67,86){${\small ba}$}
\put(103,86){${\small ba^2}$}
\put(67,145){${\small a}$}
\put(107,145){${\small a^2}$}
\put(63,5){${\small a^{-1}}$}
\put(103,5){${\small a^{-2}}$}
\put(66,65){${\small ba^{-1}}$}
\put(103,65){${\small ba^{-2}}$}
\put(212,145){${\small a^{2^{s-1}m-1}}$}
\put(160,145){${\small a^{2^{s-1}m-2}}$}
\put(212,5){${\small a^{2^{s-1}m+1}}$}
\put(160,5){${\small a^{2^{s-1}m+2}}$}
\put(160,86){${\small ba^{2^{s-1}m-2}}$}
\put(224,100){${\small ba^{2^{s-1}m-1}}$}
\put(160,65){${\small ba^{2^{s-1}m+2}}$}
\put(224,50){${\small ba^{2^{s-1}m+1}}$}
\put(240,84){${\small ba^{2^{s-1}m}}$}
\put(285,77){${\small a^{2^{s-1}m}}$}
\end{picture}
\end{center}
\caption{$\Sigma=\Cay(\mathrm{D}_{2n},S), S=\{a,a^{-1},b\}.$}
\label{figure2}
\end{figure}

Clearly, $\Sigma$ is the ladder graph of order $2n$. Let $A_1$  be the stabilizer of $1$ in $A$ and let $A_1^*$ be the subgroup of $A_1$ fixing the neighbours of $1$ in $\Sigma$, that is, $S$ pointwise.
Note that $b$ lies on two $4$-cycles of $\Sigma$ passing through $1$, but $a$ and $a^{-1}$ lie on exactly one $4$-cycle passing through $1$, respectively. This, together with the connectedness of $\Sigma$ and the transitivity of $A$ on $V(\Sigma)$, implies that $A_1^*=1$ and $|A_1|=2$.
Let $\a$ be the automorphism of $\D2n$ induced by $a\mapsto a^{-1}$ and $b\mapsto b$. Then $\a\in \Aut(\D2n,S)$ and $\Aut(\D2n,S)=\langle \a\rangle$. Since $|A_1|=2$ and $\Aut(\D2n,S)\leq A_1$, we have $A_1=\Aut(\D2n,S)$, and by Proposition~\ref{N_AUT}, $\Sigma$ is a normal Cayley graph. Furthermore, $A=R(\D2n)\rtimes \Aut(\D2n,S)$.

To prove that $\Sigma$ is not a CI-graph, by Proposition~\ref{CI-graph-prop} we only need to show that $A$ has a regular dihedral subgroup, which is not $R(\D2n)$.

Note that $R(a)^{R(b)}=R(a^{-1})$, $R(b)^\a=R(b)$ and $R(a)^\a=R(a^{-1})$. Then $R(b)\a$ is an involution and $R(a)^{(R(b)\a)}=R(a^{-1})^\a=R(a)$, that is, $R(a)$ commutes with $R(b)\a$. Since $R(a)$ has order $n$ and $n$ is even, $R(ab)\a=R(a)(R(b)\a)$ has order $n$. Furthermore, $$(R(ab)\a)^{R(b)}=R(ab)^{R(b)}\a=R(ba)\a=\a R(ba)^\a=\a R(ba^{-1})=(R(ab)\a)^{-1}.$$ Thus, $\langle R(ab)\a, R(b)\rangle$ is a dihedral group of order $2n$. If $R(\D2n)=\langle R(ab)\a, R(b)\rangle$, then $\a\in R(\D2n)$, which is impossible. Thus, $R(\D2n)\not=\langle R(ab)\a, R(b)\rangle$. To finish the proof, it suffices to show that $\langle R(ab)\a, R(b)\rangle$ is regular on $ \D2n$.

Note that $(R(ab)\a)^2=R(ab)\a R(ab)\a=R(ab)R(ab)^\a=R(ab)R(a^{-1}b)=R(a^2)$. Clearly, $R(a^2)$ has order $n/2$, and since $n\geq 3$, $\langle R(a^2)\rangle \unlhd \langle R(ab)\a, R(b)\rangle$.
Since $\langle R(a)\rangle$ is semiregular on $\D2n$ with two orbits, $\langle R(a^2)\rangle$ has four orbits on $\D2n$, that is, $\langle a^2\rangle$, $a\langle a^2\rangle$, $b\langle a^2\rangle$ and $ba\langle a^2\rangle$. The involution $R(b)$ interchanges $\langle a^2\rangle$ and  $b\langle a^2\rangle$, and $a\langle a^2\rangle$ and $ba\langle a^2\rangle$.
Furthermore, $R(ab)\a$ interchanges $\langle a^2\rangle$ and $ba\langle a^2\rangle$, and $a\langle a^2\rangle$ and $b\langle a^2\rangle$. It follows that  $\langle R(ab)\a, R(b)\rangle$ is transitive on $\D2n$, and hence regualr as $|\langle R(ab)\a, R(b)\rangle|=2n$.  \qed

\begin{cor}\rm\label{lem-S-C}
Let $\mathrm{D}_8=\langle a,b\ |\ a^4=b^2=1, bab=a^{-1}\rangle$ and $S=\{a,a^{-1},b\}$. Then  $\Cay(\mathrm{D}_8,S)$ is a non-normal non-CI-graph. In particular, $\mathrm{D}_8$ is a non-CI-group.
\end{cor}

\pf \pf Let $\Sigma=\Cay(\mathrm{D}_8,S)$ and $R=R(\mathrm{D}_8)\rtimes \Aut(\mathrm{D}_8,S)$. Then $\Aut(\mathrm{D}_8,S)=\langle \a\rangle$, where $\a$ is induced by $a\mapsto a^{-1}$ and $b\mapsto b$. By the proof of Lemma~\ref{lem-S}, $R$ contains two regular subgroups $R(\mathrm{D}_8)$ and  $\langle R(ab)\a, R(b)\rangle$, which are not conjugate in $R$.

It is easy to see that  $\Sigma\cong K_{4,4}-4K_2$, that is, the complete bipartite graph $K_{4,4}$ minus one factor. Then $\Aut(\Sigma)\cong S_4\times \mathrm{C}_2$, where $S_4$ is the symmetric group of degree $4$. Since $\Aut(\Sigma)\not=R$, $\Sigma$ is non-normal. Note that $R$ is a Sylow $2$-subgroup of $\Aut(\Sigma)$. Let $\b\in \Aut(\Sigma)$ with $o(\b)=3$. By the Sylow Theorem, all Sylow $2$-subgroups of $\Aut(\Sigma)$ are conjugate, and hence $\Aut(\Sigma)$ has exactly three Sylow $2$-subgroups: $R$, $R^\b$ and $R^{\b^2}$. Furthermore, $\Aut(\Sigma)=\langle \b\rangle R=\langle \b\rangle R^\b=\langle \b\rangle R^{\b^2}$, implying that the conjuacy class of  $R(\mathrm{D}_8)$ in $\Aut(\Sigma)$ contains three subgroups, which are subgroups of $R$, $R^\b$ and $R^{\b^2}$, respectively. Since $\Aut(\Sigma)\cong S_4\times \mathrm{C}_2$, any two distinct Sylow $2$-subgroups have intersection isomorphism to $\mathrm{C}_2^3$, so $\langle R(ab)\a, R(b)\rangle$ is not a subgroup of $R^\b$ or $R^{\b^2}$. Thus, $\langle R(ab)\a, R(b)\rangle$ is not in the conjuacy class of  $R(\mathrm{D}_8)$, that is,  $R(\mathrm{D}_8)$ and $\langle R(ab)\a, R(b)\rangle$ are not conjugate in $\Aut(\Sigma)$. By the Babai criterion in \cite{Babai}, $\Sigma$ is a non-CI-graph. \qed

\medskip

Now we are ready to prove Theorem~\ref{mainth}.

\medskip

\noindent{\bf Proof of Theorem~\ref{mainth}:} Let $n\geq 2$. First we prove (1) and (3) are equivalent, that is, $\mathrm{D}_{2n}$ is a $\rm{NDCI}$-group if and only if either $n=2,4$ or $n$ is odd. The necessity follows from Lemma~\ref{lem-S}.

To prove the sufficiency, let $n$ be odd, or $n=2,4$, and we only need to prove that $\mathrm{D}_{2n}$ is a NDCI-group. Let $\Gamma=\Cay(\mathrm{D}_{2n},S)$ be  a normal Cayley digraph. It suffices to show that $\Gamma$ is a CI-digraph. Note that we use the notations or formulae in Eqs~(\ref{eq1})-(\ref{eq6}). Let $A=\Aut(\Gamma)$. By Proposition~\ref{N_AUT}, $A=R(\mathrm{D}_{2n})\Aut(\mathrm{D}_{2n},S)\leq \Hol(\mathrm{D}_{2n})$ and $A_1=\Aut(\mathrm{D}_{2n},S)\leq \Aut(\D2n)$.

Assume $n=2$. Then $\mathrm{D}_4\cong C_2\times C_2$, $|V(\Ga)|=4$ and $A\leq S_4$. Since the symmetric group $S_4$ has a unique subgroup isomorphism to $\ZZ_2\times\ZZ_2$, $R(\mathrm{D}_{4})$ is the unique subgroup of $A$ isomorphism to $\ZZ_2\times\ZZ_2$. By the Babai criterion in \cite{Babai}, $\mathrm{D}_4$ is a DCI-group, and hence a NDCI-group. Assume $n=4$.
By Corollary~\ref{lem-S-C}, $\mathrm{D}_8$ is a non-DCI-group. However, with the help of {\sc Magma}~\cite{magma}, one may easily check that $\mathrm{D}_8$ is a NDCI-group: it also can be proved by restricting the valency of $\Ga$ to be no more than $3$ because
$\Ga\cong \Cay(\mathrm{D}_8,\mathrm{D}_8\backslash S\cup\{1\})$ and by using the fact $|A|$ is a divisor of $|\Hol(\mathrm{D}_8)|=64$.

Assume that $n$ is odd. Without loss of any generality, we may further assume that $p_1>p_2>\ldots>p_m$, where $p_i$'s are the all prime factors of $n$. For $1\leq i\leq m$, set $\pi_i=\{p_1,p_2,\ldots,p_i\}$. Recall that $A=R(\mathrm{D}_{2n})\Aut(\mathrm{D}_{2n},S)\leq \Hol(\mathrm{D}_{2n})$ and $A_1=\Aut(\mathrm{D}_{2n},S)\leq \Aut(\D2n)$. Since $C_n$ is characteristic in $\D2n$, $A_1$ fixes $C_n$ setwise. By Lemma~\ref{AutD2n}~(1), $\langle\theta_a\rangle $ is the kernel of $\Aut(\D2n)$ on $C_n$, and therefore, the kernel of $A_1$ on $C_n$ is $\langle \theta_a\rangle\cap A_1$. Then $A_1/(\langle
\theta_a\rangle\cap A_1)$ induces a subgroup of $\Aut(C_n)$, say $B$. Note that $\Aut(C_n)$ is viewed as a subgroup of $\Aut(\D2n)$, and so is $B$, too. It follows that
$$ A=R(\D2n)\rtimes A_1, \ A_1=(\langle\theta_a\rangle\cap A_1)\rtimes B \mbox{ with } B\leq \Aut(C_n).$$

Let $G$ be a regular subgroup of $A$ such that $G\cong \mathrm{D}_{2n}$. To prove that $\Gamma$ is a CI-digraph, by Proposition~\ref{CI-graph-prop} it suffices to show that $G=R(\mathrm{D}_{2n})$. We argue by contradiction, and we suppose that $G\not=R(\mathrm{D}_{2n})$.

By Lemma~\ref{holomorph} (2), $$\langle R(b)\rangle\times\Aut(C_n)_2=\langle R(b)\rangle\times\Aut(C_{p_1^{r_1}})_2\times \Aut(C_{p_2^{r_2}})_2\times\ldots\times \Aut(C_{p_m^{r_m}})_2$$ is a Sylow $2$-subgroup of $\Hol(\D2n)$. For short, we use the following notation to denote the above Sylow $2$-subgroup of $\Hol(\D2n)$:
$$HD_2=\langle R(b)\rangle\times \Aut(C_n)_2=\langle R(b)\rangle\times\Aut(C_{p_1^{r_1}})_2\times \Aut(C_{p_2^{r_2}})_2\times\ldots\times \Aut(C_{p_m^{r_m}})_2.$$

Let $B_2$ be a Sylow $2$-subgroup of $B$. Since $\Aut(C_n)$ is abelian, we have $B_2\leq \Aut(C_n)_2$, and hence $\langle R(b)\rangle\times B_2$ is a Sylow $2$-subgroup of $A$. Since Sylow $2$-subgroups of $A$  are conjugate and $|G_2|=2$,  by the Sylow theorem,  there is $d\in A$ such that $G\cap (\langle R(b)\rangle\times B_2)^d\not=1$. Thus, $G^{d^{-1}}\cap (\langle R(b)\rangle\times B_2)\not=1$. Write $H=G^{d^{-1}}$. Then $H\leq A$ is  regular on $V(\Gamma)$, and since $R(\D2n)\unlhd A$ and $R(\D2n)\not=G$, we have $$H\not=R(\D2n)\ \mbox{ and }\ |H\cap HD_2|=2.$$
Note that $H\leq A\leq \Hol(\D2n)$. By Eq~(\ref{eq6}) and Lemma~\ref{holomorph}~(1), $$\Hol(\D2n)=R(\D2n)\rtimes (\langle \theta_{a}\rangle\rtimes\Aut(C_n))=((\langle R(a)\rangle\times\langle\theta_a\rangle)\rtimes\langle R(b)\rangle)\rtimes\Aut(C_n).$$
Since $H\cong\D2n$, we may assume that
$$H=\langle v,w\rangle\cong \D2n, o(v)=n, o(w)=2, v^w=v^{-1}, \mbox { and } w\in HD_2.$$

\medskip
\noindent {\bf Claim 1:} $v\in \langle R(a)\rangle\times \langle \theta_a\rangle$.

Since $o(v)=n$, we have $o(v^{n/p_i^{r_i}})=p_i^{r_i}$. Write $T_0=\{1\}$, $T_1=\{1,v^{n/p_1^{r_1}}\}$,
$T_2=\{1,v^{n/p_1^{r_1}},v^{n/p_2^{r_2}}\}$, $\ldots$, $T_m=\{1, v^{n/p_1^{r_1}},v^{n/p_2^{r_2}},\ldots,v^{n/p_m^{r_m}}\}$. Clearly, $\langle v\rangle=\langle T_m\rangle$. To finish the proof of Claim~1, it suffices to show that $T_m\subseteq \langle R(a)\rangle\times \langle \theta_a\rangle$. To do this, we proceed by induction on $k$ to show $T_k\subseteq \langle R(a)\rangle\times \langle \theta_a\rangle$, where $0\leq k\leq m$.

Clearly, $T_0\subseteq \langle R(a)\rangle\times \langle \theta_a\rangle$, and we may let $k>0$. By induction hypothesis,
we may assume that $T_j\subseteq \langle R(a)\rangle\times \langle \theta_a\rangle$ for all $0\leq j< k$ and aim to show $T_k\subseteq \langle R(a)\rangle\times \langle \theta_a\rangle$.
Since $T_k=T_{k-1}\cup \{v^{n/p_k^{r_k}}\}$, we only need to show that $v^{n/p_k^{r_k}}\in \langle R(a)\rangle\times \langle \theta_a\rangle$.

Since $o(v^{n/p_k^{r_k}})=p_k^{r_k}$, $v^{n/p_k^{r_k}}$ is a $\pi_k$-element. By Lemma~\ref{holomorph}~(2), $v^{n/p_k^{r_k}}\in (\langle R(a)\rangle\times\langle\theta_a\rangle)\rtimes (\Aut(C_{p_1^{r_1}})_{\pi_k}\times\ldots\times \Aut(C_{p_k^{r_k}})_{\pi_k})$. Then we may write $v^{n/p_k^{r_k}}=x\b_1\b_2\ldots\b_k$, where $x\in \langle R(a)\rangle\times \langle \theta_a\rangle$ and $\b_j\in \Aut(C_{p_j^{r_j}})_{\pi_j}$ for $1\leq j\leq k$.

Clearly, $v$ commutes with every element in $T_{k-1}$, and so does $v^{n/p_k^{r_k}}$.
Since $\langle R(a)\rangle\times \langle \theta_a\rangle$ is abelian, $x$ commutes with every element in $T_{k-1}$ because $T_{k-1}\subseteq \langle R(a)\rangle\times \langle \theta_a\rangle$.  For every $1\leq \ell\leq k-1$, by Lemma~\ref{holomorph}~(3) we have that if $j\not=\ell$ then $\b_j\in \Aut(C_{p_j^{r_j}})_{\pi_j}$ commutes with every element of order $p_\ell^{r_\ell}$ in $T_{k-1}$, and then $\b_\ell$ commutes with every element of order $p_\ell^{r_\ell}$ in $T_{k-1}$  because $\b_\ell=\b_{\ell-1}^{-1}\ldots \b_1^{-1}x^{-1}v^{n/p_k^{r_k}}\b_k^{-1}\ldots \b_{\ell+1}^{-1}$, which implies that $\b_\ell$ commutes with an element of order $p_\ell^{r_\ell}$ in $\langle R(a)\rangle\times\langle \theta_a\rangle$ as $T_{k-1}\subseteq\langle R(a)\rangle\times\langle \theta_a\rangle$. Again by Lemma~\ref{holomorph}~(3), we obtain $\b_\ell=1$. It follows $v^{n/p_k^{r_k}}=x\b_k=R(y)\theta_z\b_k$ for some $y,z\in \langle a\rangle$, and since $R(y)\in A$, we have $\theta_z\b_k\in A$, so $\theta_z\b_k\in A_1$. Since $\Gamma$ is normal, Lemma~\ref{non-normal} implies that $\b_k=1$, and therefore $v^{n/p_k^{r_k}}=R(y)\theta_z\in \langle R(a)\rangle\times \langle \theta_a\rangle$, implying $T_k\subseteq \langle R(a)\rangle\times \langle \theta_a\rangle$. This completes the proof of Claim~1. \qed

By Claim~1, we may assume that $v=R(a)^k\theta_a^\ell$ for some $k,\ell\in\ZZ_n$. Since $\theta_a$ fixes $\langle a\rangle$ pointwise, we have $1^{\langle v\rangle}=\langle a^k\rangle$, the orbit of $\langle v\rangle$ containing $1$ in $\D2n$, and since $\langle v\rangle\leq H$ is semiregular, $|1^{\langle v\rangle}|=|\langle a^k\rangle|=o(v)=n$, forcing $o(a^k)=n$. Thus, $\langle a^k\rangle=\langle a\rangle$ and hence $(k,n)=1$. Since $v\in H$, if necessary we replace $v$ by a power of $v$, and then one may let $$v=R(a)\theta_a^\ell \ \ \mbox{ for some }\ell\in\ZZ_n.$$

Recall that $w\in H$ with $o(w)=2$ and $w\in HD_2=\langle R(b)\rangle\times \Aut(C_n)_2$. If $w\in \Aut(C_n)_2=\Aut(C_{p_1^{r_1}})_2\times \Aut(C_{p_2^{r_2}})_2\times\ldots\times \Aut(C_{p_m^{r_m}})_2$, then $w$ fixes $\langle a\rangle$ setwise, forcing that $\langle a\rangle$ is an orbit of length $n$ of $H=\langle v,w\rangle$,  contradicting the regularity of $H$. Thus, $w\in R(b)\Aut(C_n)_2$, that is, $$w=R(b)\e_1\e_2\ldots\e_m, \ \mbox{ where } \e_i\in \Aut(C_{p_i^{r_i}})_2\ \mbox{ and } \e_i^2=1 \mbox{ for every }\  1\leq i\leq m.$$

\medskip
\noindent {\bf Claim 2:} For every $1\leq k\leq m$, either $p_k^{r_k}\ | \ell$ and $\e_k=1$, or $(\ell,p_k)=1$ and $\e_k\not=1$.

Write $\e=\e_1\ldots\e_m\in\Aut(C_n)$. Since $H\cong\D2n$, we have $v^{R(b)\e}=v^w=v^{-1}=R(a^{-1})\theta_{a^{-\ell}}$. By Lemma~\ref{holomorph}, $$\theta_a^{R(b)}=\theta_a R(b)^{\theta_a} R(b)=\theta_a R(b^{\theta_a} b)=\theta_a R(a^{-1}) \mbox{ and } (\theta_a^\ell)^{R(b)}=\theta_a^\ell R(a^{-\ell})=R(a^{-\ell})\theta_{a^\ell}.$$
Since $v^{R(b)}=(v^{-1})^\e=R((a^\e)^{-1})\theta_{(a^\e)^{-\ell}}$, we have
$$R(a^{-(\ell+1)})\theta_{a^\ell}=R(a^{-1})R(a^{-\ell})\theta_{a^\ell}=R(a)^{R(b)}(\theta_a^\ell)^{R(b)}=v^{R(b)}=
R((a^\e)^{-1})\theta_{(a^\e)^{-\ell}}.$$ Since $\langle R(a)\rangle\cap \langle\theta_a\rangle=1$, we deduce $R(a^{\ell+1})=R(a^\e)$ and $\theta_{a^\ell}=\theta_{(a^\e)^{-\ell}}$. It follows that $a^\e=a^{\ell+1}$ and $(a^\e)^\ell=a^{-\ell}$. This yields
$a^{\ell(\ell+2)}=1$, so $n\ \mid \ell(\ell+2)$ and $p_k^{r_k}\ \mid \ell(\ell+2)$ for every $1\leq k\leq m$.
If $p_k\mid \ell$ and $p_k\mid \ell+2$, then $p_k\mid 2$, which is impossible because $n$ is odd. It follows that
either $p_k^{r_k}\mid \ell$ and $(p_k,\ell+2)=1$, or $p_k^{r_k}\mid (\ell+2)$ and $(p_k,\ell)=1$.

Assume $p_k^{r_k}\mid \ell$ and $(p_k,\ell+2)=1$. Since $a^\e=a^{\ell+1}$, we have $$a^\e=(a_1a_2\ldots a_m)^\e=a_1^{\e_1}a_2^{\e_2}\ldots a_m^{\e_m}=a^{\ell+1}=a_1^{\ell+1}\cdot a_{k-1}^{\ell+1}a_k a_{k+1}^{\ell+1}\ldots a_m^{\ell+1},$$ implying $a_k^{\e_k}=a_k$, that is, $\e_k=1$.

Assume $p_k^{r_k}\mid (\ell+2)$ and $(p_k,\ell)=1$. Since $(a^\e)^\ell=a^{-\ell}$, we have $(a_k^\ell)^{\e_k}=(a_k^\ell)^{-1}$, and so $\e_k\not=1$ because $(p_k,\ell)=1$ implies  $o(a_k^\ell)=p_k^{r_k}$. This completes the proof of Claim~2.\qed

If $p_k^{r_k}\ | \ell$ for all $1\leq k\leq m$, by Claim~2 we have $\e_k=1$, and hence $v=R(a)$ and $w=R(b)$, yielding $H=R(\D2n)$, a contradiction. Thus, $\e_k\not=1$ for some $k$. To simplify notation, from now on we do not assume $p_1>p_2>\ldots>p_m$ any more, which will not cause confusion. Then we may assume that
there exists $1\leq s\leq m$ such that $(\ell,p_1\ldots p_s)=1$, $\e_j\not=1$ for all $1\leq j\leq s$, $p_{s+1}^{r_{s+1}}\ldots p_m^{r_m}\mid \ell$,  and $\e_j=1$ for all $s+1\leq j\leq m$. It follows that $o(a^\ell)=p_1^{r_1}\ldots p_s^{r_s}$, forcing $\langle a^\ell\rangle=\langle a_1a_2\ldots a_s\rangle$ and $\langle \theta_a^\ell\rangle=\langle \theta_{a_1a_2\ldots a_s}\rangle$. Write  $\g=\e_1\ldots\e_s$. Then $w=R(b)\g$ and $\g$ is the automorphism of $\D2n$ induced by  $$a_i^\g=a_i^{-1} \mbox{ for } 1\leq i\leq s, \ \ a_j^\g=a_j \mbox{ for } s+1\leq j\leq m, \ \ b^\g=b.$$

Write $L=\langle a_1a_2\ldots\a_s\rangle$. Then $\g$ fixes every coset of $L$ in $\D2n$. Note that   $v=R(a)\theta_a^\ell\in H\leq A$ and $w=R(b)\g\in A$. Since $R(\D2n)\leq A$, we have $\theta_a^\ell\in A_1$ and $\g\in A_1$. In particular, $\langle \theta_{a_1a_2\ldots a_s}\rangle=\langle \theta_a^\ell\rangle\leq A_1$. Note that $\langle \theta_{a_1a_2\ldots a_s}\rangle$ fixes $\langle a\rangle$ pointwise, and fixes every coset of $L$ in $b\langle a\rangle$ setwise but transitive on the coset. Since $R(b)$ interchanges $\langle a\rangle$ and $b\langle a\rangle$, $\langle\theta_{a_1a_2\ldots a_s}\rangle^{R(b)}$ fixes $b\langle a\rangle$ pointwise, and fixes every coset of $L$ in $\langle a\rangle$ setwise but transitive on the coset. Since $\langle \theta_{a_1a_2\ldots a_s}\rangle\leq A$ and $\langle\theta_{a_1a_2\ldots a_s}\rangle^{R(b)}\leq A$, both the induced bipartite digraphs $[dL, cL]$ and $[cL, dL]$, for all $cL\subset\langle a\rangle$ and $dL\subset b\langle a\rangle$, are isomorphic to either the empty graph with $2p_1^{r_1}\ldots p_s^{r_s}$ isolated vertices, or the complete bipartite digraph $\vec{K}_{p_1^{r_1}\ldots p_s^{r_s},p_1^{r_1}\ldots p_s^{r_s}}$ of order $2p_1^{r_1}\ldots p_s^{r_s}$. This implies that every automorphism of the induced sub-digraph $[\langle a \rangle]$ of $\langle a\rangle$ in $\Gamma$ that fixes every coset of $L$ in $\langle a\rangle$ setwise, can be extended to an automorphism of $\Gamma$ such that it fixes $b\langle a\rangle$ pointwise.

Let $\overline{\g}$ be the permutation on $\D2n$ such that $x^{\overline{\g}}=x^\g$ for $x\in \langle a\rangle$ and $y^{\overline{\g}}=y$ for $y\in b\langle a\rangle$. Since $s\geq 1$, we have $\overline{\g}\not=1$, and since $\g\in A$, $\overline{\g}$ is an automorphism of the induced sub-digraph $[\langle a\rangle]$ fixing every coset of $L$ in $\langle a\rangle$ setwise. Then the above paragraph implies $\overline{\g}\in A$. Since $\overline{\g}$ fixes $1$, we have $\overline{\g}\in A_1=\Aut(\D2n,S)\leq \Aut(\D2n)$, and since $\overline{\g}$ fixes $b\langle a\rangle$ pointwise, we have $\overline{\g}=1$, a contradiction. This completes the proof of the
the equivalence of (1) and (3).

\medskip

Now we are ready to finish the proof of  Theorem~\ref{mainth} by proving that (2) and (3) are equivalent, that is, $\mathrm{D}_{2n}$ is a $\rm{NCI}$-group if and only if either $n=2,4$ or $n$ is odd. The sufficiency follows from the fact that a $\rm{NDCI}$-group is a $\rm{NCI}$-group, and the necessity follows from Lemma~\ref{lem-S}.
\qed

\medskip

\noindent{\bf Proof of Corollary~\ref{cor-DCI}:} Let the dihedral group $\mathrm{D}_{2n}$ be a DCI-group. Then $\mathrm{D}_{2n}$ is a NDCI-group. By Theorem~\ref{mainth}, $n$ is $2,4$ or odd. By Corollary~\ref{lem-S-C}, $\mathrm{D}_8$ is a non-CI-group and so a non-DCI-group. If $n$ is odd, by \cite[Theorem~1.2]{DE2} or \cite[Theorem~1.2]{Li2}, $n$ is square-free. Thus, $n=2$ or $n$ is odd-square-free. Now let $\mathrm{D}_{2n}$ be a CI-group. Then $\mathrm{D}_{2n}$ is a NCI-group. Again by Theorem~\ref{mainth}, either $n$ is $2,4$, or $n$ is odd. By \cite[Theorem~1.2]{DE2} or \cite[Theorem~1.2]{Li2}, if $n$ is odd, then $n=9$ or $n$ is square-free. Since $\mathrm{D}_8$ is a non-CI-group, we have $n=2, 9$, or $n$ is odd-square-free. \qed
\medskip

\medskip
\noindent {\bf Acknowledgements:} The work was supported by the National Natural Science Foundation of China (11731002, 12071023) and the 111 Project of China (B16002).

\end{document}